\documentclass{ifacconf}

\usepackage{graphicx}      
\usepackage{natbib}        
\usepackage{amsmath}
\usepackage{amssymb}
\usepackage{latexsym}
\usepackage{mathrsfs}
\usepackage{amscd}
\usepackage{sidecap}
\usepackage{graphicx}
\usepackage{color}
\usepackage{stmaryrd}
\usepackage{algorithmic}
\usepackage{algorithm}
\graphicspath{{pics/}}
\def\Q{ {\mathcal {Q}}}

\def\U{ {\mathcal {U}}}

\newcommand{\R}{\mathbb R}
\newcommand{\J}{\mathcal J}
\newcommand{\N}{\mathbb N}
\newcommand{\dt}{\Delta t}

\newcommand{\epsi}{{\varepsilon}}

\def\be#1\ee{\begin{equation}#1\end{equation}}

\newcommand{\bx}{\mathbf{x}}

\newcommand{\by}{\mathbf{y}}

\renewcommand{\theequation}{\arabic{section}.\arabic{equation}}
\renewcommand{\thetable}{\arabic{section}.\arabic{table}}
\renewcommand{\thefigure}{\arabic{section}.\arabic{figure}}
\newcommand{\bq}{\begin{equation}}
\newcommand{\eq}{\end{equation}}

\def\bqa{\begin{eqnarray}}
\def\eqa{\end{eqnarray}}
\def\dt{\Delta t}

\def\S{\sigma}

\def\mf{\mu^F}
\def\ml{\mu^L}
\newcommand{\FL}{\mathcal{K}^{FL}}
\newcommand{\FF}{\mathcal{K}^{FF}}
\newcommand{\LL}{\mathcal{K}^{LL}}

\newcommand{\bd}{\begin{displaymath}}
\newcommand{\ed}{\end{displaymath}}
\newcommand{\ba}{\begin{eqnarray}}
\newcommand{\ea}{\end{eqnarray}}

\def\L{\mathcal L}
\def\R{\mathbb{R}}
\def\N{\mathbb{N}}
\def\R{\mathbb{R}}

\def\pa{\partial}

\def\epsi{ \varepsilon}

\newcommand{\lt}{\left}
\newcommand{\rt}{\right}

\newcommand{\F}{\mathcal{F}}

\usepackage{graphicx}      
\usepackage{natbib}        
\usepackage{amsmath}
\begin{document}
\begin{frontmatter}

\title{(Sub)Optimal feedback control of mean field multi-population dynamics} 


\author[First]{Giacomo Albi} 
\author[Second]{Dante Kalise}

\address[First]{Department of Informatics, University of Verona, Str. le Grazie 15, Verona,  I-37134, Italy (e-mail: giacomo.albi@univr.it).}
\address[Second]{Department of Mathematics, Imperial College London, South Kensington Campus, London SW7 2AZ, United Kingdom (e-mail: dkaliseb@imperial.ac.uk)}

\begin{abstract}                
	We study a multiscale approach for the control of agent-based, two-population models. The control variable acts over one population of leaders, which influence the population of followers via the coupling generated by their interaction. We cast a quadratic optimal control problem for the large-scale microscale model, which is approximated via a Boltzmann approach. By sampling solutions of the optimal control problem associated to binary two-population dynamics, we generate sub-optimal control laws for the kinetic limit of the multi-population model. We present numerical experiments related to opinion dynamics assessing the performance of the proposed control design.
\end{abstract}


\end{frontmatter}
\section{Introduction and problem description} Agent-based models (ABM) constitute an adequate framework for addressing different modelling challenges in social dynamics, human crowd motion, and collective behaviour (\cite{CDFSTB03,cristiani2014book}).  While self-organization and pattern emergence phenomena has been extensively studied in the context of ABM \cite{2007TAC,motsch2011JSP},  the design of control mechanisms enforcing a certain dynamic equilibria is a novel research topic. In particular, ABM are inherently large-scale and therefore its control and optimisation is computationally demanding. For ABM where interactions are  homogeneous, a natural way to circumvent the curse of dimensionality is by means of mean field modelling, i.e., transitioning from a microscale description of the system's state to the analysis of its macro properties governed by partial differential equations (PDEs), \cite{CCR11,di2013measure,cafotove10}.  The subject of the present paper is to exploit these ideas  in order to synthesize control laws for multi-population dynamics.

We consider a two-population model, where a leaders' population ($L$) wants to influence a followers' population ($F$) toward a desired regime. Here, the state of the agent corresponds to a physical state (position, velocity) as in crowd motion application, or to an abstract characteristic of the individual like the degree of affinity with a certain statement, as in opinion dynamics.  The influence of the leaders over the followers is modelled as an optimal control problem where we minimise the distance of the populations with respect to the reference regime, constrained to the dynamic evolution of the two interacting populations. The two-population model is cast as a system of two coupled nonlinear transport PDEs
\be
\begin{aligned}\label{model}
	&\pa_t \mf + \nabla \cdot \lt(\lt(\FF[\mf] + \FL[\ml]\rt)\mf \rt) = 0,\\
	&\pa_t \ml + \nabla \cdot \lt(\lt(\LL[\ml] + u\rt)\ml \rt) = 0,
\end{aligned}
\ee
where $\mf=\mf(x,t)$ is the followers' density of agents at time $t$ at state $x\in \R^d$, and $\ml=\ml(y,t)$ is the leaders' density at time $t$ with state $y\in\R^d$. We assume that these densities have different total mass, respectively
\be
\rho^F = \int_{\R^d} \mf(x,t) dx,\quad \rho^L = \int_{\R^d} \ml(y,t) dy,
\ee
which are conserved quantities in time.
The interaction forces $ \mathcal{K}^{ij}[\cdot]$ for  $i,j\in\{F,L\}$ are described by the nonlocal operator
\begin{equation}\label{eq:convolution}
\mathcal{K}^{ij}[\mu^j](x,t) = \int_{\R^d} K^{ij}(x,y)(y-x)\mu^j(y,t)\,dy\,
\end{equation}
where the binary kernel $K^{ij}(\cdot,\cdot)$ encodes social interaction rules between agents such as attraction, repulsion, or alignment. The term $u=u(y,t)$ corresponds to a dynamical strategy, which we design by minimising an energy measure of the followers and leaders states
\begin{align}\label{eq:MFJ}
u^*:=&\underset{u\in\U}{argmin}\int\limits_{0}^{T}e^{-\lambda t} \L(t,\mf,\ml,u) \ dt,
\end{align}
where $\L$ stands for a performance measure which typically enforces a reference consensus position, see for example \cite{APZa,ACFK16,BBCK18}, the control signal belongs to an admissible $\U$, and $\lambda\geq 0$ is a discount factor.
Models of type \eqref{model}, have been also derived as a mean field approximation of agent-based systems of $N$ followers and $M$ leaders. Indeed, the complexity of casting an optimal control problem for a coupled system of $N+M$ non-linear differential equations becomes computationally prohibitive when $N+M$ becomes large. On the other hand, the optimal control of system \eqref{model} requires refined numerical techniques to resolve consistently the non-local integro-differential operators, and which have to cope with the solution of the optimal control problem. The approximation of the convolution operator \eqref{eq:convolution} at every time step is costly even for low-dimensional microscopic state spaces (of dimension $d$).

In order to cope with the control of problem of \eqref{model}, we rely on the approximation of the mean field dynamics by means of binary type interactions. This approach has been developed for the control of a single mean field population dynamics in the seminal paper \cite{AHP}, combining {\em model predictive control}, \cite{camacho2010model} applied to binary interaction system and the {\em quasi-invariant scaling} of Boltzmann-type equations, \cite{T}. Extension to models with interacting populations has been studied in \cite{during2009PRSA}, and in \cite{APZa} for control dynamics. Whereas, most recently, in \cite{ACFK16,AFK17} we derived a hierarchy of optimal feedback controllers based on binary dynamics, which shows a better coherence with respect to the optimal control solution of a mean field single-population model. Hence, in the present work we want to extend this novel hierarchy of controls to the case of mean field multi-population models. 

The paper is structured as follows. In Section 2 we introduce the microscopic two-population model, its control and the optimal feedback synthesis for the binary system. In Section 3 we present the Boltzmann-Bellman approximation, which allow us to compute (sub)optimal controllers for the mean field model upon the solution of binary control problems. Finally, Section 4 is devoted to the presentation of numerical approximation method for the realisation of the proposed feedback law and the closed-loop dynamics.

\section{The microscale two-population model and its control}
We consider a two-population model of $N$ \textsl{followers} agents, each represented by $x_i(t)\in\R^d$, and $M$ \textsl{leaders} $y_j(t)\in\R^d$, evolving according to
\begin{align}
 \frac{dx_i}{dt}&= \frac{1}{N}\sum_{k=1}^{N} P(x_i,x_k)(x_k-x_i)\cr
 &\qquad+\frac{1}{M}\sum_{l=1}^{M} S(x_i,y_l)(y_l-x_i)\,,\label{eq:MAS1}\\
  \frac{dy_j}{dt}&= \frac{1}{M}\sum_{l=1}^{M} R(y_j,y_l)(y_l-y_j)+u(t)\,,\label{eq:MAS2}\\
   x_i(0)&=x_0\,,\qquad i=1,\ldots,N\,,\label{eq:MAS3}\\
 y_j(0)&=y_0\,,\qquad j=1,\ldots,M\,,\label{eq:MAS4}
\end{align}
where $P(\cdot,\cdot), R(\cdot,\cdot), S(\cdot,\cdot): \R^d\times\R^d\longrightarrow\R^d$ are interaction functions (follower-follower, follower-leader, and leader-leader, respectively). For the sake of simplicity, in this section we use the $(P,R,S)$ notation for the kernels instead of $(K^{FF},K^{FL},K^{LL})$ used in other sections . The control variable $u(t)\in\U=\{u(t): \R_+\longrightarrow U\}$, with $U$ a compact subset of $\R^d$ is a signal broadcasted to the leader population only. We denote by $\bx(t)=(x_1(t),\ldots,x_{N}(t))^\top$, and $\by(t)=(y_1(t),\ldots,y_{M}(t))^\top$.
The control signal is synthesized upon the following optimal control problem 
\begin{align}\label{eq:cost}
\underset{u(\cdot)\in\U}{\min} \J(u(\cdot);\bx_0,\by_0):=\int_0^\infty e^{-\lambda t}\ell(\bx(t),\by(t),u(t))\,dt\,,
\end{align}
with a discount factor $\lambda>0$, subject to the system dynamics \eqref{eq:MAS1}-\eqref{eq:MAS4}. The running cost $\ell(\bx,\by,u)$ is of the form
\begin{align}
\ell(\bx,\by,u):=\frac{a_F}{N}\|\bx-\bar x\|_2^2+\frac{a_L}{M}\|\by-\bar x\|_2^2+\gamma\|u\|_2^2,\label{eq:penalizations}
\end{align}
with $\gamma\geq0$, $\bar x \in\R^{N\times d}$ is a desired reference state, and
\begin{align}
\quad\|\bx\|^2_2=\sum_{i=1}^N|x_i|^2\,,\quad\|\by\|^2_2=\sum_{j=1}^M|y_j|^2\,,
\end{align}
where $|\cdot|$ stands for the $d$-dimensional Euclidean norm. For the sake of simplicity, in the following we restrict our analysis to the case $d=1$, although the presented methodology is directly applicable to multidimensional agent systems. The different penalizations in \eqref{eq:penalizations} represent the deviation of the leaders with respect to a reference value, the cohesiveness of the overall leader-follower population, and control energy, respectively. It is through the cohesiveness term that leaders act steering the followers to the desired reference.

We shall  assume that system dynamics have been discretized in time with a first-order approximation

\begin{align}
 \frac{x_i^{n+1}-x_i^n}{\dt}&=\frac{1}{N}\sum_{k=1}^{N} P^n_{ik}(x_k^n-x_i^n)+\frac{1}{M}\sum_{l=1}^{M} S^n_{il}(y_l^n-x_i^n)\,,\label{eq:MAS1dt}\\
\frac{y_j^{n+1}-y_j^{n}}{\dt}&=\frac{1}{M}\sum_{l=1}^{M} R^n_{jl}(y_l^n-y_j^n)+u^n\,,\label{eq:MAS2dt}
\end{align}
where $P^n_{ij}:=P(x^n_i,x^n_j),\, S_{i,l}^n =S(x_1^n,y_l^n)$, $R^n_{ij}:=R(y^n_i,y^n_j)$, and $x^n=x(n\dt)$, with $n\in\N$ and a time discretisation parameter $\dt>0$. The cost functional \eqref{eq:cost} is discretized accordingly
\begin{align}
\J_{\dt}(u;\bx^0,\by_0):=\sum_{n=0}^{\infty} \beta^n\ell(\bx^n,\by^n,u^n)\,,
\end{align}
with $\beta=e^{-\lambda\dt}$. For this infinite horizon control problem we shall focus on solutions of \eqref{eq:cost} which can be expressed in feedback form
\begin{align}\label{feedcontr}
u^*(t)=\F_{\dt}(\bx(t),\by(t))\,,
\end{align}
i.e. controllers which can be computed based on the present information of the system. For this, we follow a dynamic programming approach, where the feedback mapping $\F_{\dt}$ is obtained from the solution of a nonlinear Bellman equation.  The Bellman equation is solved over a domain of the same dimension of the state space, and therefore the applicability of the dynamic programming approach is limited to low-dimensional dynamics. In order to circumvent this difficulty, we resort to the analysis of control problems over binary dynamics, i.e. two-particle systems, which will be conveniently used to generate a control law for the mean field ABM.

\subsection{Optimal feedback synthesis for the binary system}
We focus our analysis on the optimal control problem when $N=M=2$. In this case, denoting by $\bx_{12} = (x_1,x_2)$ and $\by_{12}  = (y_1,y_2)$, we have the following binary control problem:
\begin{align}\label{bincost}
\underset{u\in\U}{\min}\;\;  \J_{\dt}(u;\bx_{12}^0,\by_{12}^0):=\sum_{n=0}^{\infty} \beta^n\ell(\bx_{12} ^n,\by_{12} ^n,u^n)\,
\end{align}
subject to the two-agent models
\begin{align}\label{binary}\nonumber
x_1^{n+1}&= x_1^n + \frac{\dt}{2}\left( P^n_{12}(x^n_2-x^n_1) + \sum_{l=1}^2 S_{1l}^n(y_l^n-x^n_1)\right)\,,\\\nonumber
x_2^{n+1}&= x_2^n + \frac{\dt}{2}\left( P^n_{21}(x^n_1-x^n_2) + \sum_{l=1}^2 S_{2l}^n(y_l^n-x^n_2)\right)\,,\\\nonumber
y_1^{n+1}&= y_1^n + \frac{\dt}{2} R^n_{12}(y^n_2-y^n_1) + \dt\, u^n\,,\\
y_2^{n+1}&= y_2^n + \frac{\dt}{2} R^n_{21}(y^n_1-y^n_2) + \dt\, u^n\,.
\end{align}

By defining the value function associated to the infinite horizon discrete cost \eqref{bincost} as
\begin{align}\label{eq:Val}
V(\bx_{12}^0,\by_{12}^0) := \underset{u\in\U}{\inf}\sum_{n=0}^{\infty}  \beta^n \ell(\bx_{12}^n,\by_{12}^n,u^n)\,,
\end{align}
then the application of the Dynamic Programming Principle characterizes the value function as the solution of the Bellman equation
\begin{equation}
\begin{aligned}\label{eq:DP}
V(\bx_{12},\by_{12}) & = \underset{u\in U}{\min}\left\{\beta V(\bx_{12}^+, \by_{12}^+(u))+\dt\ell(\bx_{12},\by_{12},u)\right\}\,,
\end{aligned}
\end{equation}
where $\bx_{12}^+$ and $\by_{12}^+$ denote a one-step update of system \eqref{binary} departing from $\bx_{12}$ and $\by_{12}$ with control action $u$.
The optimal feedback mapping is recovered from $V$ via the relation
\begin{equation}\label{eq:ochj}
u^*=\underset{u\in U}{argmin}\left\{\beta V(\bx_{12}^+, \by_{12}^+(u))+\dt\ell(\bx_{12},\by_{12},u)\right\}
\end{equation}

\section{The Boltzmann-Bellmann approximation}
In order to consistently scale the microscopic constrained dynamics to a large-scale MAS, we approximate the feedback systems introduced in \eqref{binary} by means of binary operators,whose evolution is described by a system of two Boltzmann equations.  Thus, we decompose system \eqref{binary} with the following binary interactions transformation,
\begin{subequations}
\label{BBfull}
\begin{align}
&
\begin{cases}\label{FF}
x_1'= x_1 + \alpha_{} K^{FF}(x_1,x_2)(x_2-x_1)\,\\
x_2'= x_2 + \alpha_{} K^{FF}(x_2,x_1)(x_1-x_2)
\end{cases}
\\&
\begin{cases}\label{FL}
x_1''= x_1 + \alpha_{} K^{FL}(x_1,y_k)(y_k-x_1)\,\\
y_k''= y_k, \quad k\in{\{1,2\}}
\end{cases}
\\&
\begin{cases}\label{LL}
y_1'= y_1 + \alpha_{} K^{LL}(y_1,y_2)(y_2-y_1) + 2\alpha_{}\, \phi_\alpha(y_1,y_2)\,,\\
y_2'= y_2 + \alpha_{} K^{LL}(y_2,y_1)(y_1-y_2) + 2\alpha_{}\, \phi_\alpha(y_2,y_1)\,,
\end{cases}
\end{align}
\end{subequations}
where $x_1',x_2',y_1'y_2'$  are the post-interactions states, and $\alpha>0$ is the analogous of the time discretisation step, $\Delta t/2$, which we will refer as interaction strength. The control $\phi_\alpha(\cdot)$ is computed averaging the feedback control \eqref{feedcontr} with respect to the followers' distribution, 
\be\label{binctrl}
\begin{aligned}
&\phi_\alpha(y_1,y_2;\mf)  =\cr
 &\qquad\iint_{\R^{2d}}\mathcal{F}_\alpha(x_1,x_2,y_1,y_2)\mf(x_1)\mf(x_2)\ dx_1 dx_2.
\end{aligned}
\ee
Hence, in order to describe at a larger scale the evolution of a system of agents' populations interacting through binary rules of type \eqref{BBfull}, we introduce the following two coupled Boltzmann equations
\begin{align}\label{eq:BoltzFL}
\begin{aligned}
&\pa_t \mu^F = \mathcal{Q}^{FF}(\mu^F,\mu^F) + \mathcal{Q}^{FL}(\mu^F,\mu^L) ,
\\
&\pa_t \mu^L =  \mathcal{Q}^{LL}(\mu^L,\mu^L),
\end{aligned}
\end{align}
with initial data $\mu^F(x,0) = \mu^F_0(x), \mu^L(x,0) = \mu^L_0(x)$, and where each quadratic operator $\mathcal{Q}^{ij}(\mu^{i},\mu^{j})$ contains the information of the interaction among populations, with $i,j\in\{F,L\}$. More precisely, $\mathcal{Q}^{ij}(\mu^{i}\mu^{j})$ is defined as follows 
\begin{align}\label{eq:Qij}
\begin{aligned}
\mathcal{Q}^{ij}(\mu^i,\mu^j) = \int_{\R^{d}}  \eta_{ij}\left(\frac{1}{\mathbb{J}_{ij}}\mu^{j}('y)\mu^{i}('x)- \mu^{j}(y)\mu^{i}(x)\right)dy
\end{aligned}
\end{align}
with $('x,'y)$ the pre-interaction state that generates state $(x,y)$ according to the interaction \eqref{BBfull}, where $\mathbb{J}_{ij}$ denotes the Jacobian of the transformation $('x,'y)\to(x,y)$ and $\eta_{ij}>0$ the interaction rate. The operator $\mathcal{Q}^{ij}$ is also named {\em loss-gain} term, since it accounts the change-rate of states. For further details on models of type \eqref{eq:BoltzFL} in the context of consensus dynamics, we refer to \cite{T,during2009PRSA,APZa}.

Finally, we  show that system \eqref{eq:BoltzFL} is equivalent to a mean field system of type \eqref{model}, under a {\em quasi invariant scaling}, \cite{T}. 
\begin{thm}\label{thm:grazing}
Let $\alpha\geq0$ and $t\geq0$, and assume $u_\alpha(\cdot)$, $K^{ij}(\cdot,\cdot)\in L^2_{loc}$, for every $i,j\in\{L,F\}$.  Furthermore we assume that for $\alpha\to 0$, the limit $\phi_\alpha(\cdot) \to \phi(\cdot)$ is well defined. Given $\epsi>0$ and the {\em quasi-invariant} scaling defined as follows
\be\label{scaling}
\alpha_{} = \varepsilon ,\quad \eta_{ij}=1/\varepsilon,\qquad 
\ee 
 for every $i,j\in\{F,L\}$, then let $(\mu^{L,\varepsilon},\mu^{F,\varepsilon})$ be a solution for the scaled system \eqref{eq:BoltzFL}, with initial datum $(\ml_0(x),\mf_0(x))$.
Then in the limit $\varepsilon\to0$, the solution $\ml_\varepsilon,\mf_\varepsilon$ converges point-wise, up to a subsequence, to  $(\ml,\mf)$ satisfying the following nonlinear mean field system
\be
\begin{aligned}\label{eq:MFFL}
&\pa_t \mf + \nabla \cdot \lt(\lt(\FF[\mf] + \FL[\ml]\rt)\mf \rt) = 0,\\
&\pa_t \ml + \nabla \cdot \lt(\lt(\LL[\ml] + \Phi[\mf,\ml]\rt)\ml \rt) = 0,
\end{aligned}
\ee
 with initial datum $(\ml_0(x),\mf_0(x))$, and  where
\be
\begin{aligned}\label{Fctrl}
\Phi[\mf,\ml](x)  =\ \int_{\mathbb{R}^{d}}\phi(x,y;\mf)\ml(z)\,dy. 
\end{aligned} 
\ee
\end{thm}
We refer to \cite{ACFK16,APZa}, for a detailed proof of the Theorem.

\section{Numerical approximation and tests}
{We propose a numerical scheme for system \eqref{eq:MFFL}, which extend the approach proposed in \cite{ACFK16},  consisting of two main parts. First, an off-line procedure for the synthesis of the feedback map \eqref{feedcontr} via the solution of \eqref{eq:DP}, approximated with a semi-Lagrangian, policy iteration scheme \cite{AFK15,KKK16}. Second, the direct solution of the Boltzmann system \eqref{eq:BoltzFL} via a modified Direct Simulation Monte-Carlo method (DSMC),~\cite{BN,PTa}. More specifically, we consider the following first-order scheme for the scaled Boltzmann system \eqref{eq:BoltzFL}
\begin{equation}\label{eq:MCBoltz}
\begin{aligned}
        &\mu^{F,n+1}=\left(1-\frac{\delta t}{\varepsilon}(\rho^F+\rho^L)\right)\mu^{F,n}\cr
        &\quad +\frac{\delta t\rho^F}{\varepsilon}\frac{{\Q}_\varepsilon^{FF,+}(\mu^{F,n},\mu^{F,n})}{\rho^F}+\frac{\delta t\rho^L}{\varepsilon}\frac{\Q_\varepsilon^{FL,+}(\mu^{F,n},\mu^{L,n})}{\rho^L},
        \\
        &\mu^{L,n+1}=\left(1-\frac{\delta t\rho^L}{\varepsilon}\right)\mu^{L,n}+\frac{\delta t\rho^L}{\varepsilon}{\frac{{\Q}_\varepsilon^{LL,+}(\mu^{L,n},\mu^{L,n})}{\rho^L}},
        \end{aligned}
\end{equation}
for $n=0,\ldots,N_T-1$, with $\mu^i,n=\mu^i(x,n\delta t)$ for $i\in\{F,L\}$, and time step $\delta t >0$ such that $N_t\delta t =T$. Note also that we split the interaction operator $\mathcal{Q}^{ij}$ in a gain and a loss part according to \eqref{eq:Qij} and the scaling \eqref{scaling},
\[
\Q_\varepsilon^{ij}(\mu^i,\mu^j) = \frac{\rho^j}{\varepsilon}\left(\frac{\Q_\varepsilon^{ij,+}(\mu^{i},\mu^{j})}{\rho^j} - \mu^{i}\right), \ {i,j}\in\{F,L\},
\]
where we have assumed that the masses $\rho^{F},\rho^L$ are conserved. Since $\mu^{L,n},\mu^{F,n}$ are positive densities, we have the following time step restriction,
 \be\label{CFL}
 \delta t\leq\varepsilon/(\rho^F+\rho^L).
 \ee
In the spirit of DSMC methods, we sample $N_s$ and $M_s$ particles respectively from $\mf(\cdot,t^n)$ and $\ml(\cdot,t^n)$, and we simulate their evolution according to system \eqref{eq:MCBoltz}.

Combining both methods, we  propose the following \emph{Two Populations Boltzmann-Bellmann} algorithm (TPBB) 
\begin{alg}[\textbf{TPBB algorithm}]~\,
  \begin{enumerate}
  \item[\texttt{0.}]  The infinite horizon feedback control, $\F_\epsi(\cdot)$, is derived on the  computational domain $\Omega\subset \R^{4d}$.
  \item[\texttt{1.}]
  Given $N_s$ and $M_s$ samples $\left\{x^0_k\right\}_{k=1}^{N_s}, \left\{y^0_r\right\}_{r=1}^{M_s}$, respectively from the initial distribution $\mf_0(x),\ml_0(y)$.
  \item[\texttt{2.}]
  \texttt{for} $n=0$ \texttt{to} $N_{T}-1$\\
  \ $-$ {\em Followers}
  \begin{enumerate}
  \item[\texttt{a.}]  Set 
  $N^{FF}_c = \left\llbracket \frac{\delta t\rho^F}{\epsi}{N_s}\right\rrbracket$, and $N^{FL}_c =\left\llbracket\frac{\delta t\rho^L}{\epsi}{N_s}\right\rrbracket$;
  \item[\texttt{b.}]  select $N^{FF}_c$ random samples uniformly without repetition among all possible $N_s$ samples; 
  \item[\texttt{c.}]  for each $x_i^n$, compute the post-interaction position $x_i'$ via \eqref{FF}, selecting a random sample $x_r^n$ among the $N_s$ followers, with repetition;
  \item[\texttt{e.}]  select $N^{FL}_c$ random samples uniformly without repetition among the $N_s-N^{FF}_c$ remaining followers,
  \item[\texttt{f.}]   for each $x_j^n$ compute the post-interaction position $x_j''$ via \eqref{FL}, selecting a random sample $y_r^n$ among the $M_s$ leaders, with repetition;
 \item[\texttt{g.}]   set $x_i^{n+1}=x_i'$ and $x_j^{n+1} = x_j''$ for all $i$ and $j$, and $x_k^{n+1} = x_k^{n}$ for the remaining $N_s - N^{FF}_c-N^{FL}_c$.
  \end{enumerate}
  \ $-$ {\em Leaders}
   \begin{enumerate}
    \item[\texttt{a.}] Compute the feedback control \eqref{binctrl} using $\sigma_s$ particles, sampled from $\mu^F(x,t^n)$,
  \be\label{numctrl}
  \phi_\epsi(\cdot,\cdot) = \frac{1}{\sigma_s^2}\sum_{h=1}^{\sigma_s}\sum_{k=1}^{\sigma_s} \F_\epsi(x_h^n,x_k^n,\cdot,\cdot);
  \ee
  \item[\texttt{b.}]  set $M^{LL}_c = \left\llbracket\frac{\delta t \rho^L}{\epsi}M_s\right\rrbracket$;
 \item[\texttt{c.}]  select $M^{LL}_c$ random samples uniformly without repetition among all possible $M_s$ samples; 
  \item[\texttt{d.}]  for each $y_i^n$, compute the post-interaction position $y_i'$ via \eqref{LL}, selecting a random agent $y_r^n$ among the $M_s$ leaders, with repetition;
  \item[\texttt{e.}] set $y_i^{n+1}=y_i'$ for all $i$ and $y_j^{n+1} = y_j^{n}$ for the remaining $M_s - M^{LL}_c$ leaders.
  \end{enumerate}
\item[]\texttt{end for}
  \end{enumerate}
 \label{ANMC}
\end{alg}\medskip
Where $\llbracket\cdot\rrbracket$ indicates the integer stochastic rounding. For further details on these algorithms we refer to \cite{PTa}.

}
\paragraph*{\bf Computational tests.} In order to show the validity of the presented approach, we propose several numerical examples for different type of the interaction kernels in the context of opinion dynamics. 
In every test we consider two populations of followers and leaders, respectively with total density $\rho^F = 1$, and $\rho^L=0.5$. We consider one-dimensional dynamics defined on  $\Omega = [-1,1]$.
Followers' and leaders' densities are reconstructed using respectively $N_s=10^6$, and $M_s = 5\times10^5$ samples. The computational domain is discretized with the parameter $\delta x  =  0.025$, whereas the time step discretisation is set to be $\delta t = 2/3\times10^{-2}$ according to \eqref{CFL}, and having fixed the scaling parameter $\epsi = 0.01$. 
We consider the infinite horizon control problem associated to functional \eqref{bincost}, with running cost  \eqref{eq:penalizations} for the reduced binary dynamics. For each test we summarize in Table \ref{tab:par} the parameters associated to \eqref{eq:penalizations}.


\begin{table}[h]
\caption{Parameters' choice for functional \eqref{bincost} for the different Test 1-2-3.}
\label{tab:all_parameters}
\begin{center}
\begin{tabular}{ccccccc }
\hline
  &$T$& $a_F$ & $a_L$ & ${\gamma}$ & $\lambda$ &  $\bar x$   \\
\hline
\hline
Test 1: & 2.5& $1$ &$1$& $1$ &$1$ & $-0.5$   \\
\hline
Test 2:& 10& $10$ &$0.1$ &$0.05$ & $0.1$ & $0.25$ \\
\hline
Test 3: & 3.5 & 1 &0.01 &$1$ & $0.5$ & $0$   \\
\hline
\end{tabular}\label{tab:par}
\end{center}
\end{table}

\subsubsection{Test 1: Linear model.}
\begin{figure}[t]
\centering
\includegraphics[scale=0.215]{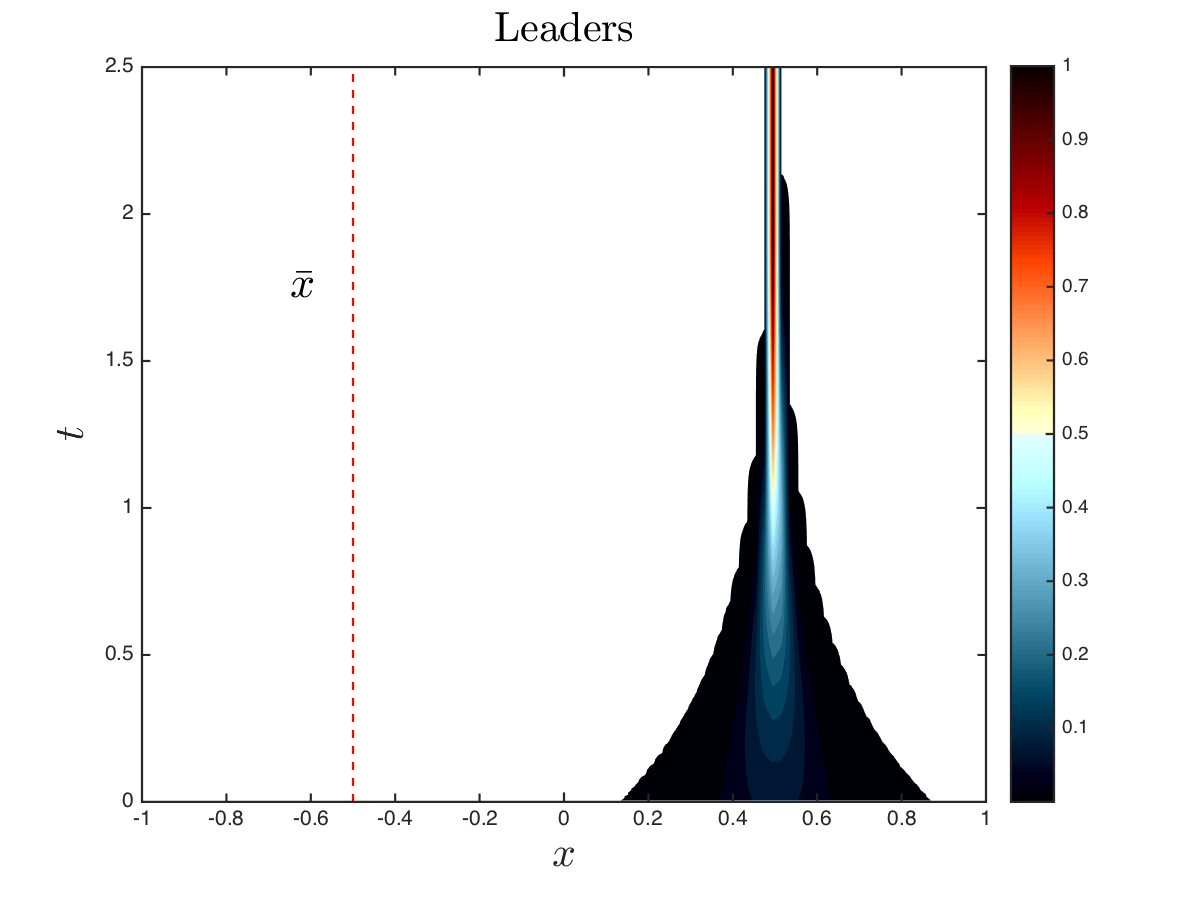}
\includegraphics[scale=0.215]{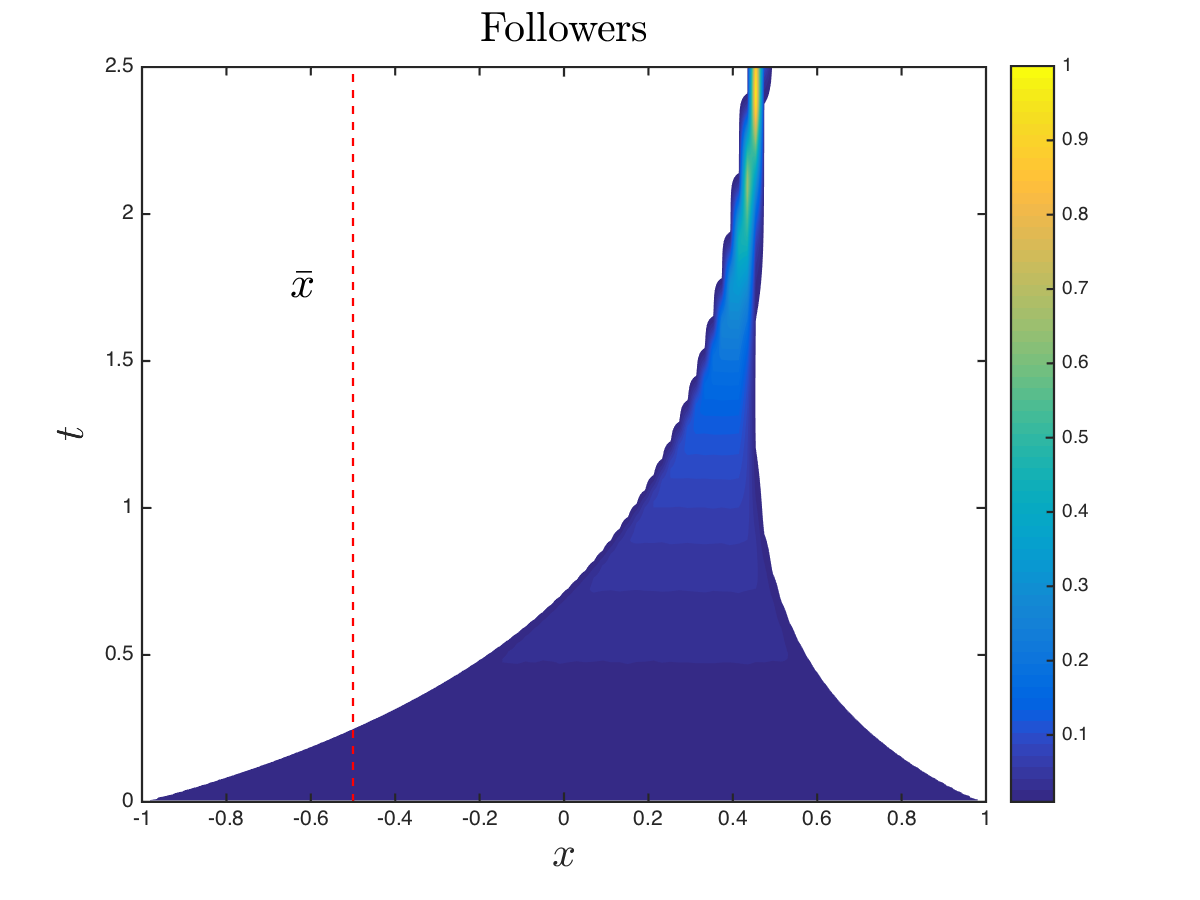}
\\
\includegraphics[scale=0.215]{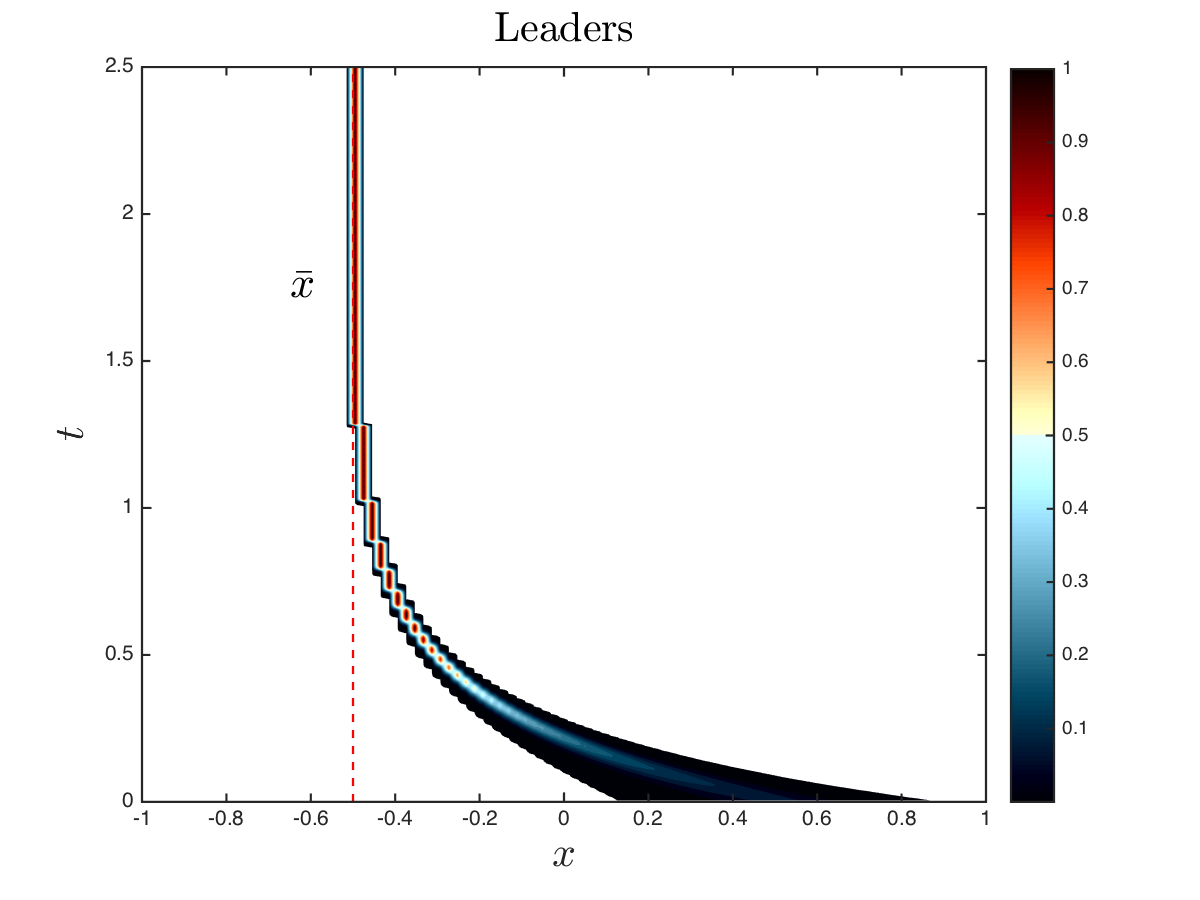}
\includegraphics[scale=0.215]{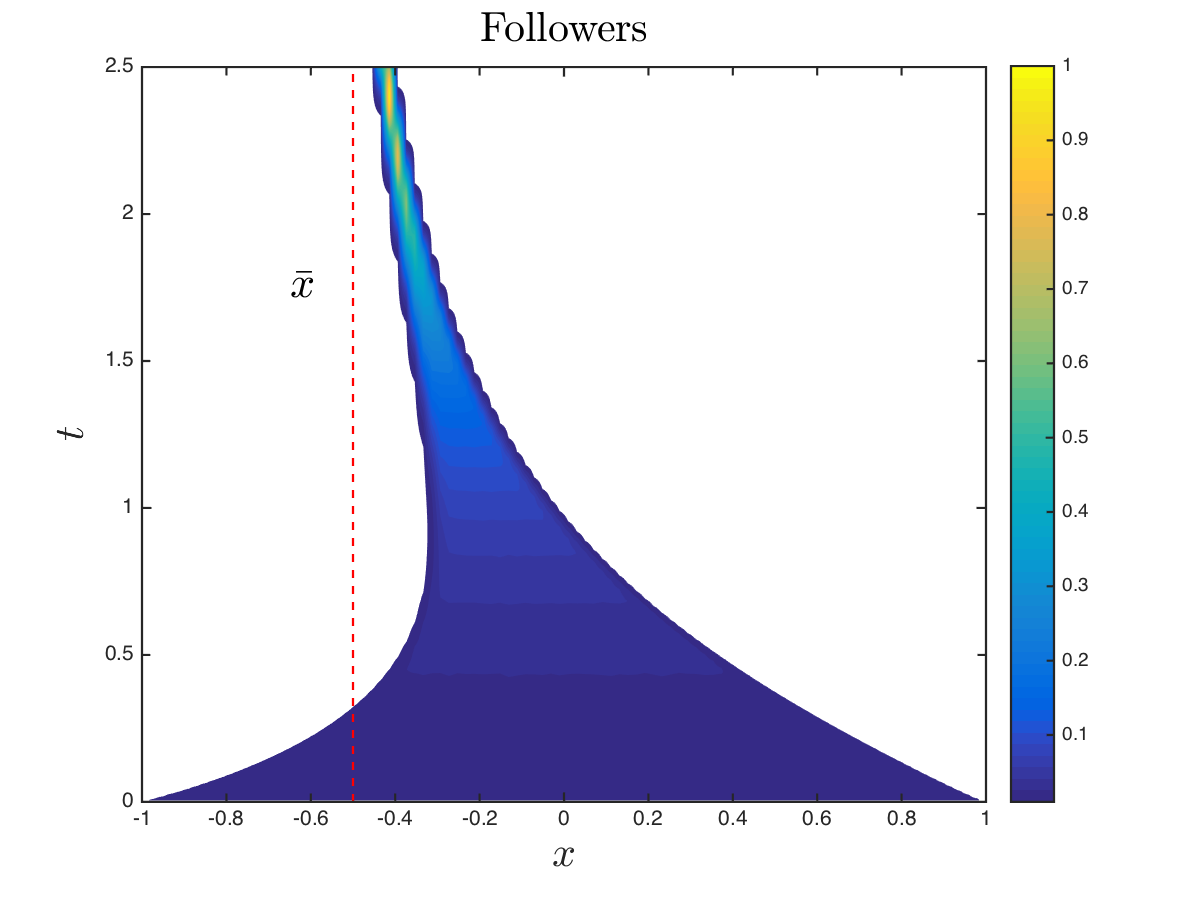}
\caption{Test 1. 
First row represents the free dynamics of the linear followers-leaders system (no control action). 
Second row shows the solution of the infinite horizon control problem for the linear model of the followers-leaders system. On the left, evolution of the leaders' density, on the right evolution of the followers' density, the control is capable to steer the system towards the reference $\bar x = -0.5$.
}
\label{Fig_T1}
\end{figure}

We consider the optimal control of the leader-follower consensus system \eqref{eq:MFFL} with linear interaction kernels, i.e. we choose the following constant interaction rate
\be
\begin{aligned}\label{eq:MFSz2}
&K^{FF}(x,y)  =  K^{FL}(x,y)  = K^{LL}(x,y)  = 1.
\end{aligned}
\ee
Note that in this case the mdel \eqref{eq:MFFL} reduces to the following linear system of coupled PDEs
\be
\begin{aligned}\label{model_L}
&\pa_t \mf + \nabla \cdot \lt(\lt((m^F(t)-x) + (m^L(t)-x)\rt)\mf \rt) = 0,\\
&\pa_t \ml + \nabla \cdot \lt(\lt(m^L(t)-x) + \Phi\rt)\ml \rt) = 0,
\end{aligned}
\ee
where $m^F,m^L$ are respectively the followers' and leaders' average opinions (state). Hence, in this linear-quadratic  case the optimal feedback map for the binary dynamics can be solved explicitly via a Riccati equation. We fix the final time $T = 2.5$, and we consider an initial uniform distribution for the followers' population  on the entire domain, $\mu^F_0(x) \sim \textrm{Unif}([-1,1])$, whereas for the leaders' population we set a subset of the opinions' domain as follows, $\mu^L_0(x) \sim \textrm{Unif}([0.15,0.85])$. In the presence of a control $u$ the leaders aim to steer the system towards the reference position $\bar x = -0.5$, and we set the control to be bounded into the set  $U=[-1,1]$. We depict the evolution of the followers' and leaders' densities in Figure \ref{Fig_T1}: the first row shows the unconstrained dynamics, whereas the second row the constrained evolution with control parameters $\gamma=1,\lambda=1$,$a_L=a_F=1$. In both cases followers are always attracted to the leaders' position. However, in the uncontrolled setting the leaders align around their initial average state $m_0^L$,  whereas in the controlled case they move towards the reference $\bar x$. 

\subsubsection{Test 2: Bounded confidence interactions.}
\begin{figure}[t]
\centering
\includegraphics[scale=0.215]{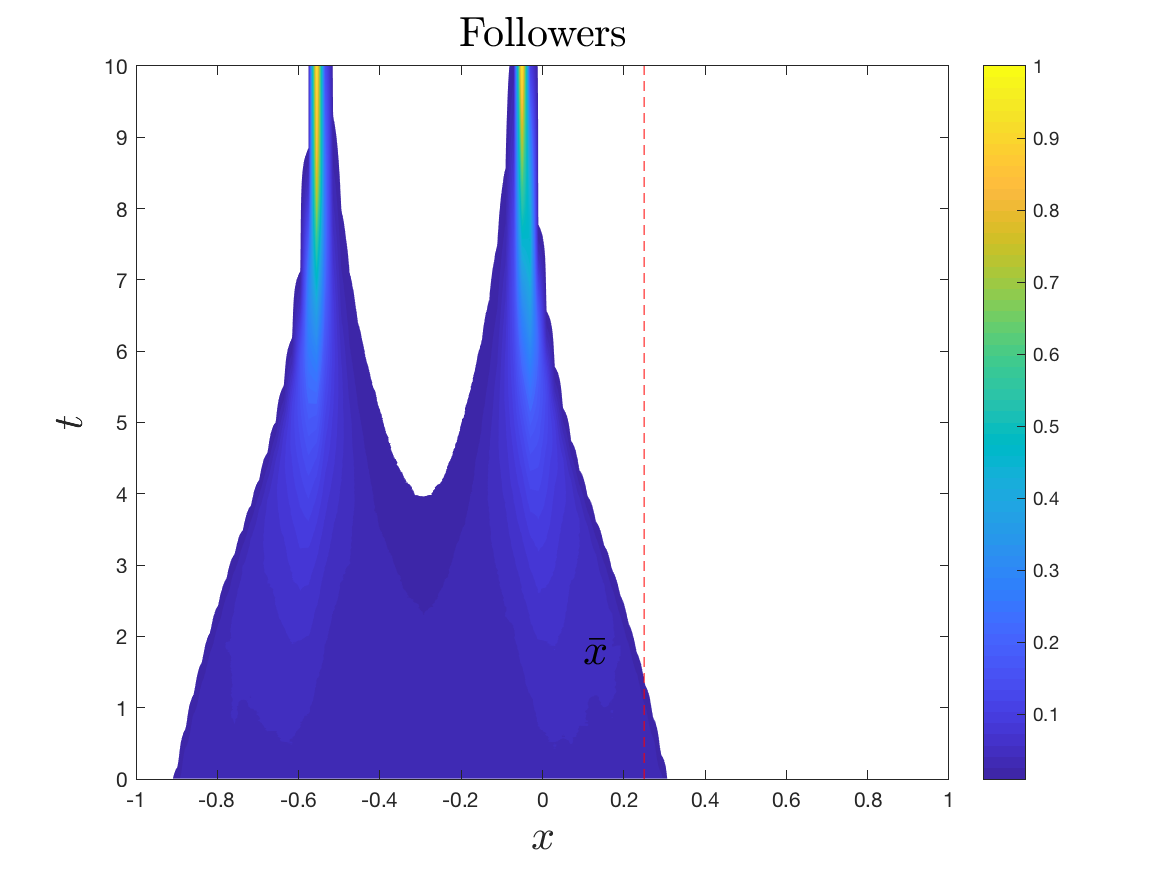}
\includegraphics[scale=0.215]{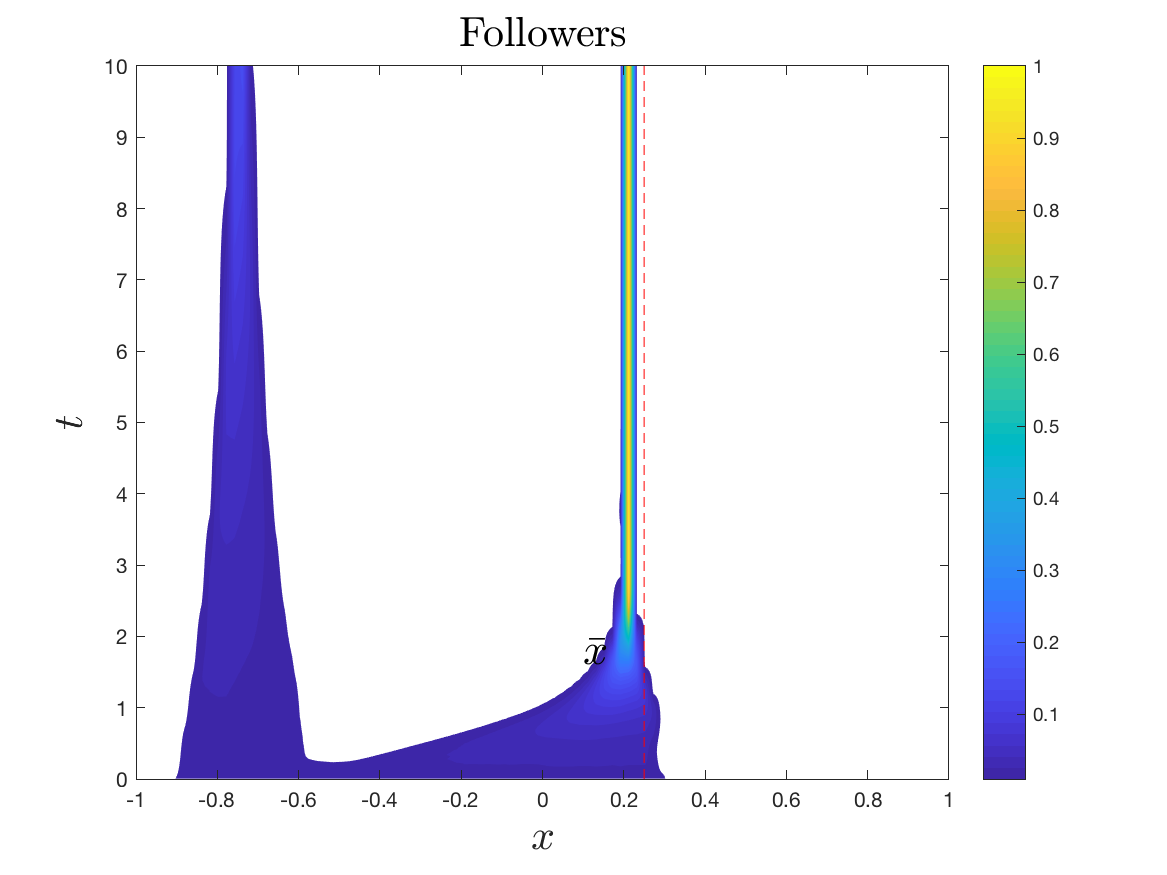}
\caption{Test 2. 
Left-hand figure represents the evolution of the unconstrained followers' density, without leaders' interaction, i.e. $K^{FL}\equiv 0$, on the right-hand side we observe how the leaders are able to steer large part of the followers towards $\bar x = 0.5$.
}
\label{Fig_T2a}
\end{figure}

We consider a non-linear interaction for followers' dynamics, while the leaders act as an aggregator of the opinions through their constrained evolution. Hence, we consider the following interactions functions: for followers we have
\be
\begin{aligned}\label{eq:FTest2}
&K^{FF}(x,y) =  \chi_{\left\{|x-y|\leq 0.3\right\}}(y),\cr  &K^{FL}(x,y)  = \chi_{\left\{|x-y|\leq 0.8\right\}}(y),
\end{aligned}
\ee
and for leaders we prescribe a linear interaction kernel $K^{LL} = 1$. Note that the followers' kernels account for the interactions among followers and followers-leaders only below a certain level of confidence. This type of interaction is also known as {\em bounded confidence model}, \cite{hekr02}.
We fix the final time $T = 10$, and we consider an initial uniform distribution for followers' population  as follows, $\mu^F_0(x) \sim \textrm{Unif}([-0.9,1.3])$, whereas leaders' population are distributed uniformly as, $\mu^L_0(x) \sim \textrm{Unif}([0,0.5])$. In the presence of control the leaders aim at steering the system towards the reference position $\bar x = 0.25$.  Here we set the control to be bounded into the set  $U=[-1,1]$, and we choose the following parameters for the infinite horizon control problem, $\lambda = 0.1$, $\gamma = 0.05$. We penalize differently the position of the followers and leaders with respect to the reference position, we set $a_F = 10$, $a_L = 0.1$.
\begin{figure}[t]
\centering
\includegraphics[scale=0.215]{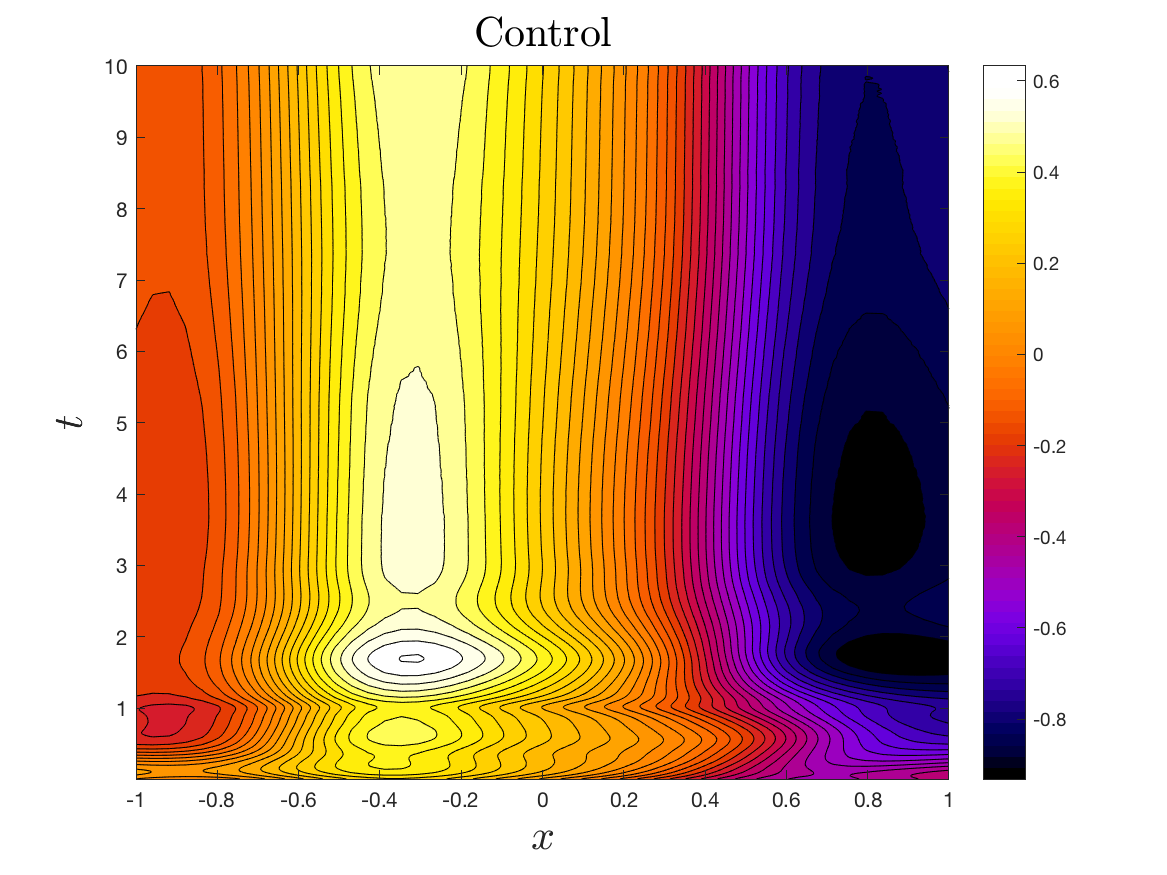}
\includegraphics[scale=0.215]{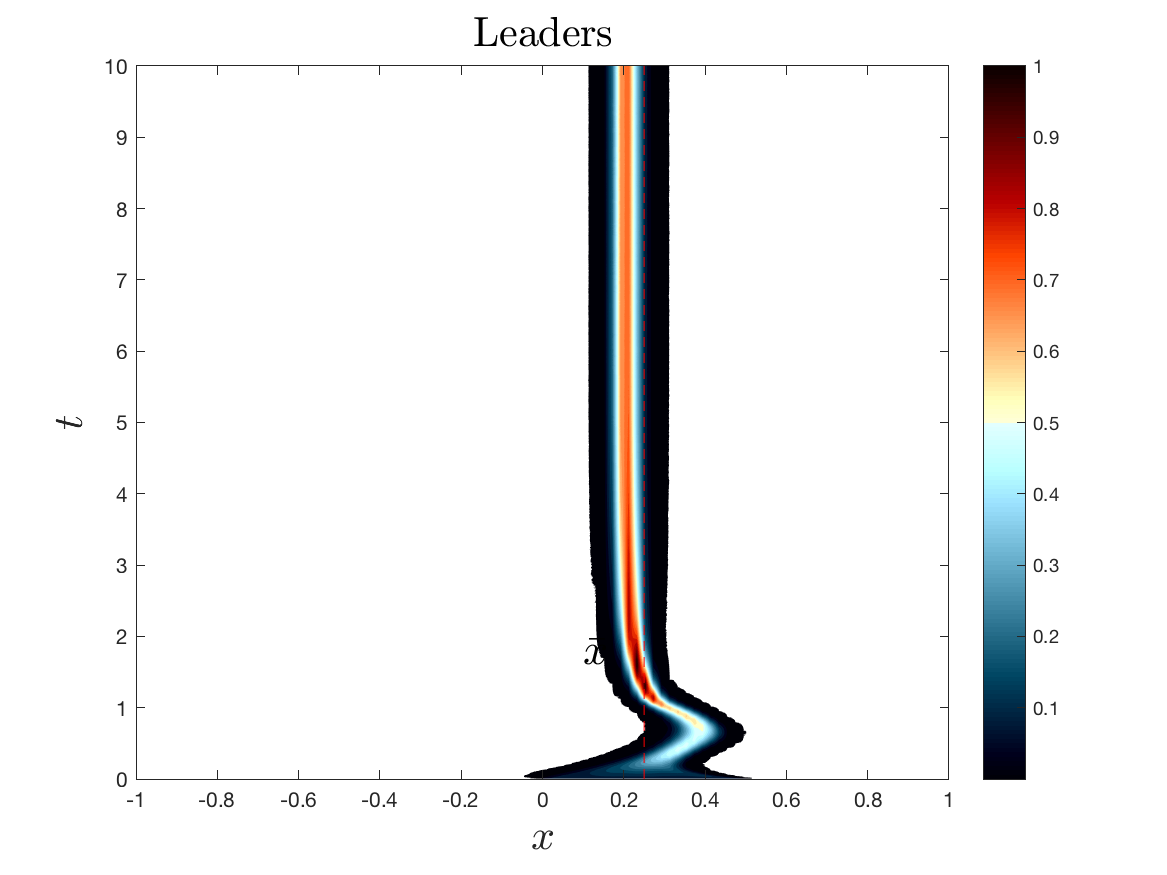}
\caption{Test 2. 
On the left-hand side we depict the shape of the control $\Phi(y,t)$.
On the right-hand side constrained evolution of the leaders' density is reported.
}
\label{Fig_T2b}
\end{figure}
On the left-hand side of Figure \ref{Fig_T2b} we report the shape of the control $\Phi(y,t)$, defined by \eqref{Fctrl}, and computed according to \eqref{numctrl} with $\sigma_s = 2\times10^5$ samples, on the right-hand side we depict the evolution of the leaders' density constrained to the control. 
Figure \ref{Fig_T2a} shows the evolution of the followers' density in two cases: on the left-hand side we have the dynamics in absence of followers-leaders interaction, i.e. $K^{FL} \equiv 0$, where we observe the emergence of two opinion's clusters; on the right-hand side, a large portion of the followers' density moves toward the reference position $\bar x = 0.25$, due to the action of the leaders.

\subsubsection{Test 3: Controlling disagreement.}
\begin{figure}[t]
\centering
\includegraphics[scale=0.215]{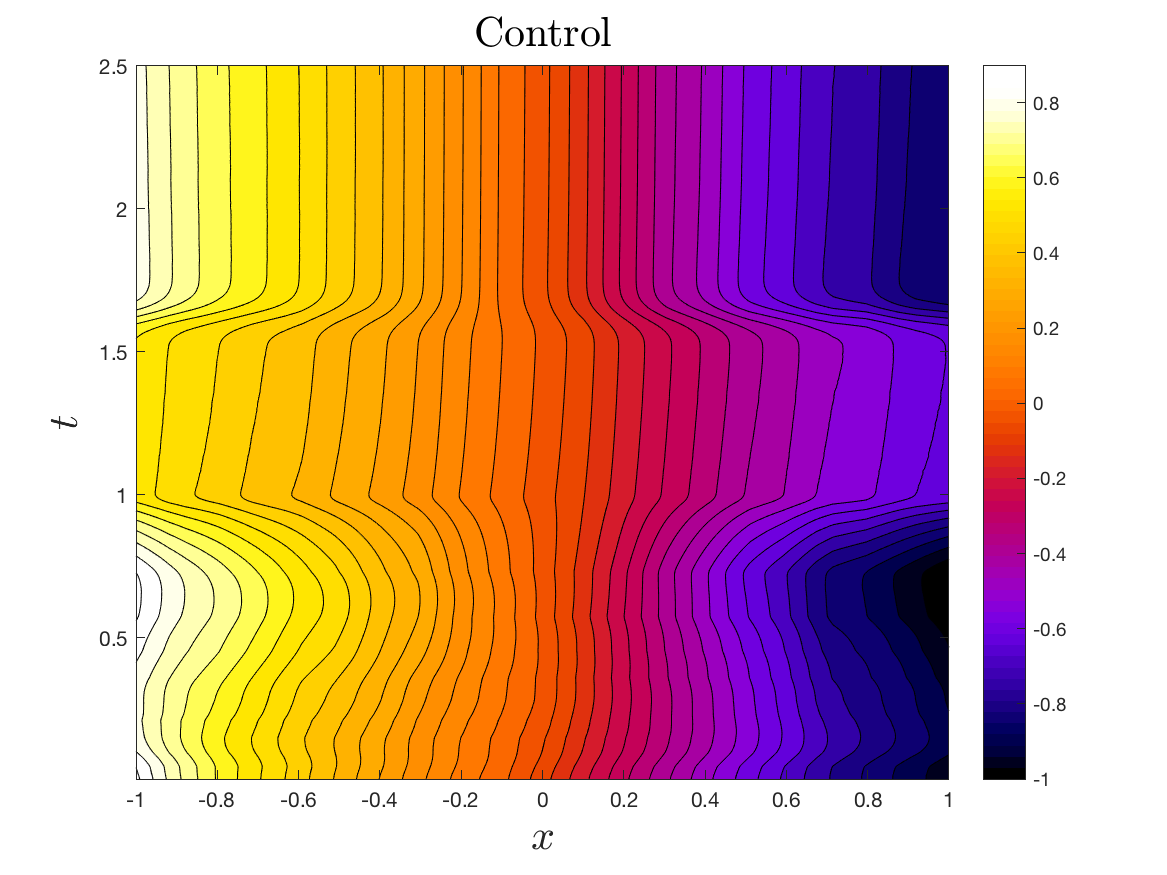}
\includegraphics[scale=0.215]{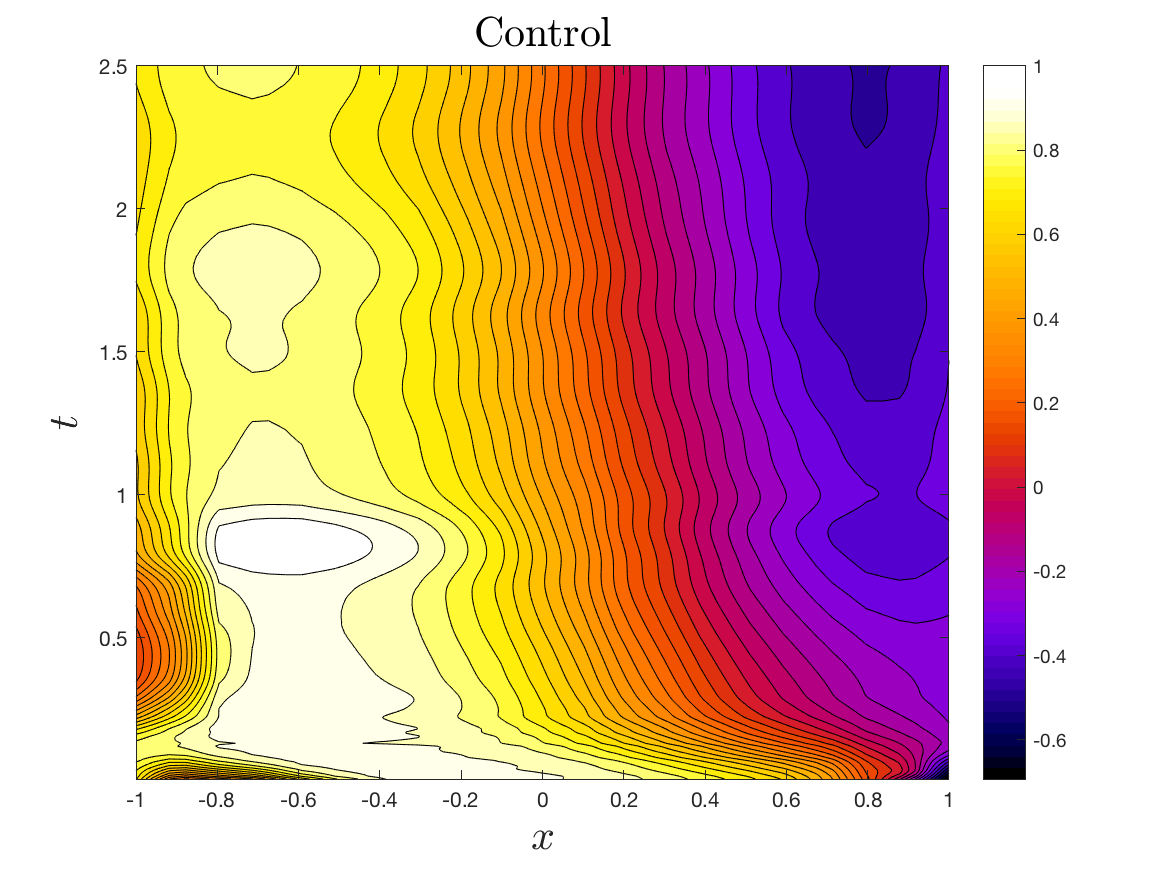}
\caption{Test 3 (Control). Shape of the control $\Phi(x,t)$ in the frame $[-1,1]\times[0,T]$. On the left-hand side the control for the attractive leaders and repulsive followers, \eqref{eq:FTest3a}, on the right-hand side the case with attractive followers and repulsive leaders.
}
\label{Fig_T31}
\end{figure}

In this final test we consider two different type of interactions for the followers' and leaders' dynamics. In the first one the followers 
are repelled by other followers, but attracted by leaders, and a second interaction where, conversely, followers are repelled by leaders and aggregate with other followers.  Hence we consider two different type of interactions for followers
\begin{subequations}
\begin{align}
&K_{}^{FF}(x,y)  = +(1-x^2), K_{}^{FL}(x,y)  = -(1-x^2),\label{eq:FTest3a}
\\
&K_{}^{FF}(x,y)  = -(1-x^2), K_{}^{FL}(x,y)  = +(1-x^2),\label{eq:FTest3b}
\end{align}
\end{subequations}
whereas leaders, in both cases, have an interaction governed by $$K^{LL}(y,z) = (1-y^2).$$
We  remark that such kernels weight the interaction with other agents according only to their own opinion, \cite{SWS}. At time $t=0$ the followers' density is distributed uniformly in the interval $I = [0.05,.55]$, and leaders are distributed uniformly as $\mu_0^L\sim\textrm{Unif}(-0.45,0.05)$. The optimal feedback control $\Phi(y,t)$ is computed using the parameters $\gamma = 1$ and $\lambda = 0.5$, whereas, similarly to the previous test we penalize less  the leaders' positions with respect to the reference position $\bar x = 0$. We set $a_F= 1$ and $a_L = 0.01$.
Having set the final time to be $T = 2.5$, in Figure \ref{Fig_T31} we report the shape of the control $\Phi(x,t)$  for the both cases reconstructed from \eqref{numctrl} using $\sigma_s = 5\times10^5$ samples: on the right-hand side for a repulsive followers and attractive leaders, \eqref{eq:FTest3a}, on the left-hand side for attractive followers and repulsive leaders, \eqref{eq:FTest3b}.

\begin{figure}[t]
\centering
\includegraphics[scale=0.215]{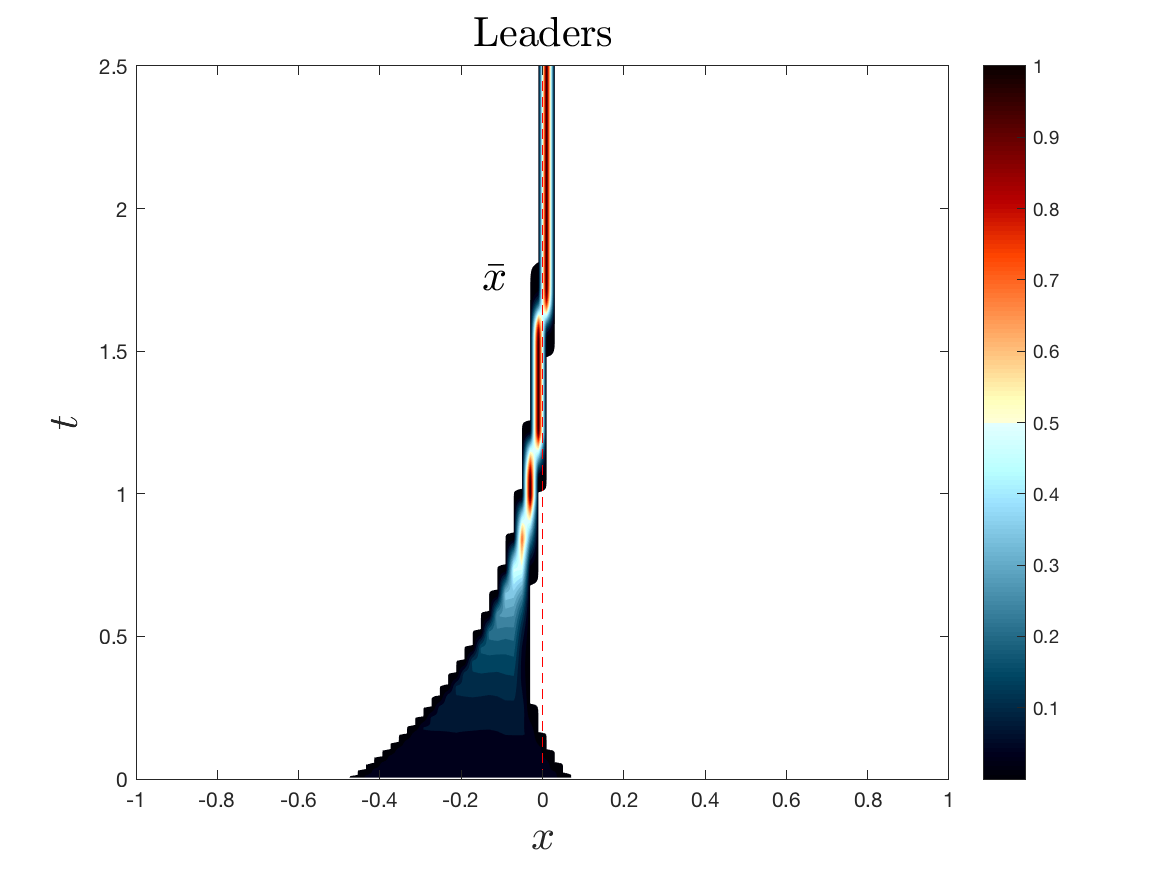}
\includegraphics[scale=0.215]{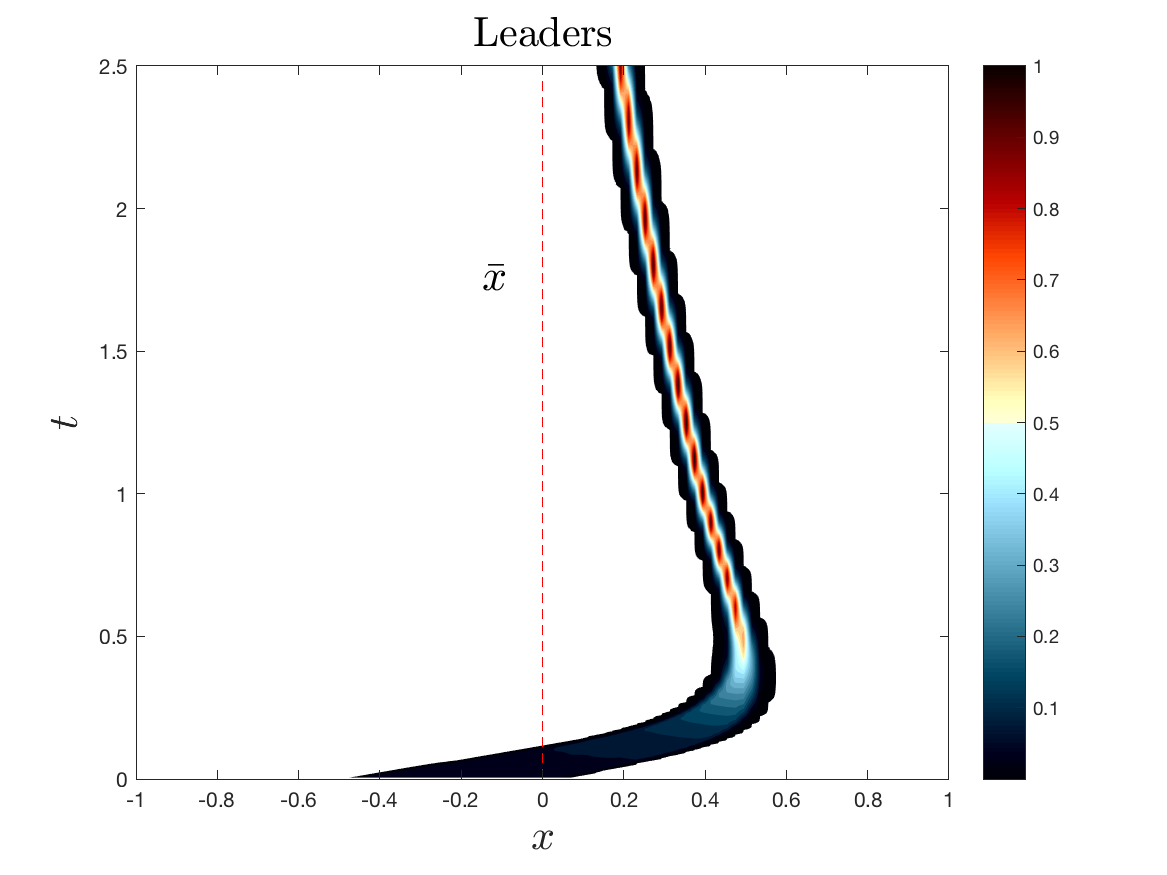}
\\
\includegraphics[scale=0.215]{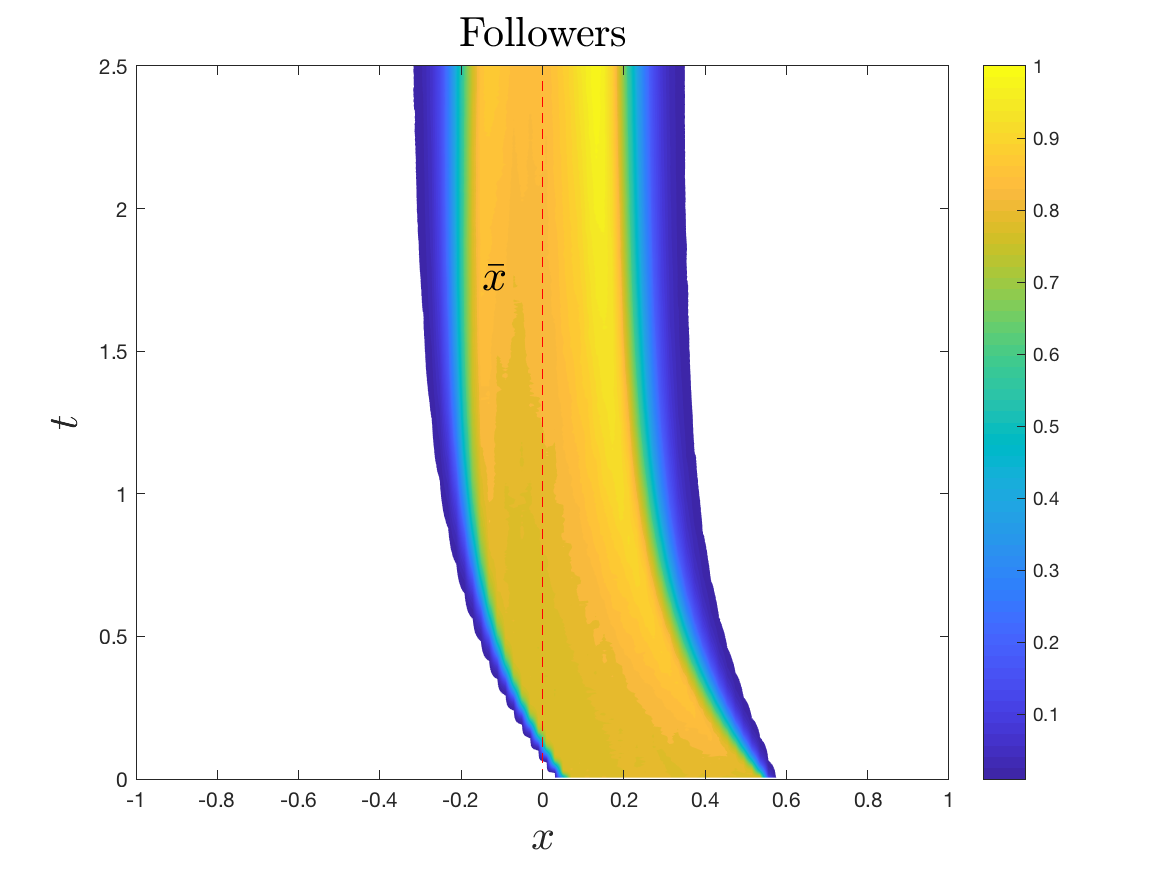}
\includegraphics[scale=0.215]{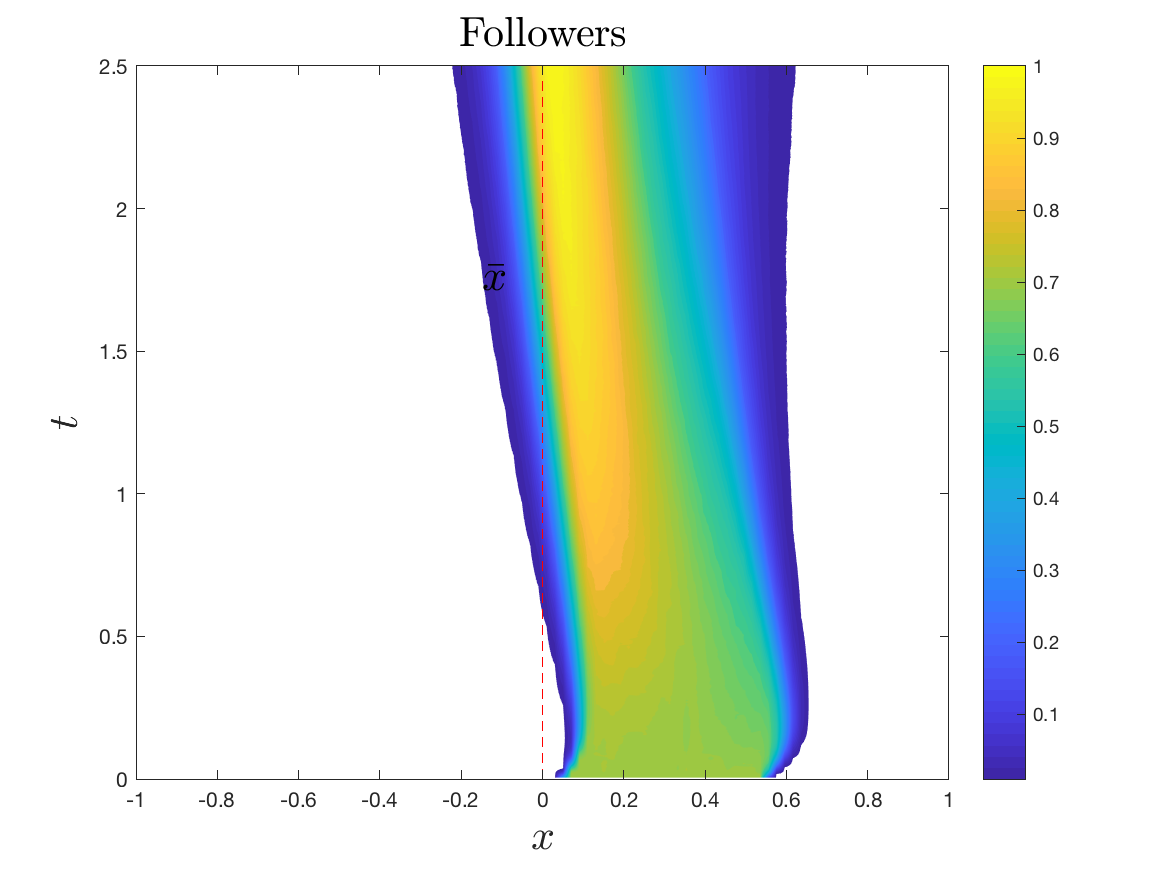}
\caption{Test 3. Evolution of the leaders-followers dynamics. On the left-hand column we report the evolution for dynamics of type\eqref{eq:FTest3a}, on the right-hand column the evolution for dynamics of type\eqref{eq:FTest3a}. 
In both cases leaders' action is capable to  influence the evolution of followers' density towards the reference position $\bar x  = 0$.
}
\label{Fig_T32}
\end{figure}

We observe in Figure \ref{Fig_T32} the evolution of the followers' and leaders' densities. The first row depict the evolution of the leaders and the second row for the followers, whereas the left-hand side figures show followers-leaders dynamics \eqref{eq:FTest3a}, and the right-hand side dynamics \eqref{eq:FTest3b}. We remark how the leaders' evolution behaves differently in order to make use of the interaction with followers. In particular, in the first case \eqref{eq:FTest3a},  the behaviour is similar to the one observed in the first test. In the second case \eqref{eq:FTest3b}, leaders initially move away from the reference position, in order to push the followers' density towards $\bar x$, thanks to their repulsive interaction.
%

\end{document}